# LEARNING BY MIRROR AVERAGING

By A. Juditsky, P. Rigollet and A. B. Tsybakov

*Université Grenoble 1, Georgia Institute of Technology and Université Paris 6 CREST and Université Paris 6*

Given a finite collection of estimators or classifiers, we study the problem of model selection type aggregation, that is, we construct a new estimator or classifier, called aggregate, which is nearly as good as the best among them with respect to a given risk criterion. We define our aggregate by a simple recursive procedure which solves an auxiliary stochastic linear programming problem related to the original nonlinear one and constitutes a special case of the mirror averaging algorithm. We show that the aggregate satisfies sharp oracle inequalities under some general assumptions. The results are applied to several problems including regression, classification and density estimation.

**1. Introduction.** Several problems in statistics and machine learning can be stated as follows: given a collection of $M$ estimators, construct a new estimator which is nearly as good as the best among them with respect to a given risk criterion. This target is called model selection (MS) type aggregation, and it can be described in terms of the following stochastic optimization problem.

Let $(\mathcal{Z}, \mathfrak{F})$ be a measurable space and let $\Theta$ be the simplex

$$\Theta = \left\{ \theta \in \mathbb{R}^M : \sum_{j=1}^{M} \theta^{(j)} = 1, \theta^{(j)} \geq 0, j = 1, \ldots, M \right\}.$$

Here and throughout the paper we suppose that $M \geq 2$ and we denote by $z^{(j)}$ the $j$th component of a vector $z \in \mathbb{R}^M$. We denote by $[z^{(j)}]_{j=1}^{M}$ the vector $z = (z^{(1)}, \ldots, z^{(M)})^\top \in \mathbb{R}^M$.

Let $Z$ be a random variable with values in $\mathcal{Z}$. The distribution of $Z$ is denoted by $P$ and the corresponding expectation by $E$. Suppose that $P$

---

is unknown and that we observe $n$ i.i.d. random variables $Z_1, \ldots, Z_n$ with values in $\mathcal{Z}$ having the same distribution as $Z$. We denote by $P_n$ the joint distribution of $Z_1, \ldots, Z_n$ and by $E_n$ the corresponding expectation.

Consider a measurable function $Q : \mathcal{Z} \times \Theta \to \mathbb{R}$ and the corresponding average risk function

$$A(\theta) = EQ(Z, \theta),$$

assuming that this expectation exists for all $\theta \in \Theta$. Stochastic optimization problems that are usually studied in this context consist in minimization of $A$ on some subsets of $\Theta$, given the sample $Z_1, \ldots, Z_n$. Note that since the distribution of $Z$ is unknown, direct (deterministic) minimization of $A$ is not possible.

For $j \in \{1, \ldots, M\}$, denote by $e_j$ the $j$th coordinate unit vector in $\mathbb{R}^M$: $e_j = (0, \ldots, 0, 1, 0, \ldots, 0) \in \mathbb{R}^M$, where 1 appears in $j$th position.

The aim of MS aggregation is to "mimic the oracle" $\min_j A(e_j)$, that is, to construct an estimator $\tilde{\theta}_n$ measurable with respect to $Z_1, \ldots, Z_n$ and called aggregate, such that

$$(1.1) \qquad E_n A(\tilde{\theta}_n) \leq \min_{1 \leq j \leq M} A(e_j) + \Delta_{n,M},$$

where $\Delta_{n,M} > 0$ is a remainder term that should be as small as possible. Thus, the stochastic optimization problem associated to MS aggregation is

$$\min_{\theta \in \{e_1, \ldots, e_M\}} A(\theta).$$

As an example, one may consider the loss function of the form $Q(z, \theta) = \ell(z, \theta^\top H)$ where $\ell : \mathcal{Z} \times \mathbb{R} \to \mathbb{R}$ and $H = (h_1, \ldots, h_M)^\top$ is a vector of preliminary estimators (classifiers) constructed from a training sample which is supposed to be frozen in our considerations (thus, $h_j$ can be viewed as fixed functions). The value $A(e_j) = E\ell(Z, h_j)$ is the risk corresponding to $h_j$. Inequality (1.1) can then be interpreted as follows: the aggregate $\tilde{\theta}_n^\top H$, that is, the convex combination of initial estimators (classifiers) $h_j$, with the vector of mixture coefficients $\tilde{\theta}_n$ measurable with respect to $Z_1, \ldots, Z_n$, is nearly as good as the best among $h_1, \ldots, h_M$. The word "nearly" here means that the value $\min_j A(e_j)$ is reproduced up to a reasonably small remainder term $\Delta_{n,M}$. Lower bounds can be established showing that, under some assumptions, the smallest possible value of $\Delta_{n,M}$ in a minimax sense has the form

$$(1.2) \qquad \Delta_{n,M} = \frac{C \log M}{n},$$

with some constant $C > 0$; cf. [24].

Besides being in themselves precise finite sample results, oracle inequalities of the type (1.1) are very useful in adaptive nonparametric estimation.

They allow one to prove that the aggregate estimator $\tilde{\theta}_n^\top H$ is adaptive in a minimax asymptotic sense (and even sharp minimax adaptive in several cases; for more discussion see, e.g., [18]).

The aim of this paper is to obtain bounds of the form (1.1)–(1.2) under some general conditions on the loss function $Q$. For two special cases [density estimation with the Kullback–Leibler (KL) loss, and regression model with squared loss] such bounds have been proved earlier in the works of Catoni [7, 8, 9] and Yang [29]. They independently obtained the bound for density estimation with the KL loss, and Catoni [8, 9] solved the problem for the regression model with squared loss. Bunea and Nobel [5] improved the regression with squared loss result of [8, 9] in the case of bounded response, and obtained some related inequalities under weaker conditions. For a problem which is different but close to ours (MS aggregation in the Gaussian white noise model with squared loss) Nemirovski [18], page 226, established an inequality similar to (1.1), with a suboptimal remainder term. Leung and Barron [15] improved upon this result to achieve the optimal remainder term.

Several other works provided less precise bounds than (1.1)–(1.2), with $K \min_j A(e_j)$ where the leading constant $K > 1$, instead of $\min_j A(e_j)$ in (1.1) and with a remainder term which is sometimes larger than the optimal one (1.2); a detailed account can be found in the survey [4] or in the lecture notes [17]. We mention here only some recent work where aggregation of arbitrary estimators is considered: [1, 6, 16, 22, 28, 30]. These results are useful for statistical applications, especially if the leading constant $K$ is close to 1. However, the inequalities with $K > 1$ do not provide valid bounds for the excess risk $E_n A(\tilde{\theta}_n) - \min_j A(e_j)$, that is, they do not show that $\tilde{\theta}_n$ approximately solves the stochastic optimization problem.

Below we study the mirror averaging MS aggregate $\hat{\theta}_n$ which is defined by a simple recursive procedure (cf. Section 3). This procedure outputs a convex mixture of initial estimators. Before defining the procedure, we give some arguments in favor of considering mixtures rather than selectors. Selectors are estimators with values in $\{e_1, \ldots, e_M\}$, for example, minimizers of the empirical risk. In Proposition 2.1 we show that selectors cannot satisfy (1.1)–(1.2), even for the simplest case where the loss function $Q$ is quadratic. The main results of the paper are given in Section 4; there we prove that the suggested mirror averaging aggregate satisfies oracle inequalities (1.1)–(1.2) under some general assumptions on $Q$. Finally, we show in Section 5 that these assumptions are fulfilled for several statistical models including regression, classification and density estimation.

**2. Suboptimality of selectors.** Recall that our goal is to construct an estimator $\tilde{\theta}_n$ that satisfies an oracle inequality of the type (1.1). A traditional way to approach this problem is based on empirical risk minimization. Define the empirical risk $A_n$ by

$$A_n(\theta) = \frac{1}{n} \sum_{i=1}^{n} Q(Z_i, \theta)$$

and the empirical risk minimizer (ERM) by

$$\tilde{\theta}_n^{\mathrm{ERM}} = \operatorname*{arg\,min}_{\theta \in \{e_1, \ldots, e_M\}} A_n(\theta).$$

Clearly, the ERM selects one of the $M$ initial estimators. More generally we call *selector* any estimator $T_n$ based on the sample $(Z_1, \ldots, Z_n)$ having this property, that is, such that $T_n$ takes values in $\{e_1, \ldots, e_M\}$.

The following example shows that under the squared loss the rate of convergence $\Delta_{n,M}$ in (1.1) for any selector $\tilde{\theta}_n = T_n$ is not faster than $\sqrt{(\log M)/n}$ which is substantially worse than the optimal rate given in (1.2).

Indeed, consider the squared loss

$$(2.1) \qquad Q(z, \theta) = \tfrac{1}{2}\theta^\top \theta - z^\top \theta, \qquad z \in \mathbb{R}^M, \ \theta \in \Theta.$$

For $k = 1, \ldots, M$ denote by $P^k$ the distribution of a Gaussian random vector $Z \in \mathbb{R}^M$ with mean $e_k(\sigma/2)\sqrt{(\log M)/n}$ and the covariance matrix $\sigma^2 I$ where $I$ stands for the identity matrix, and denote by $E^k$ the corresponding expectation. It is easy to see that the risk $A_k(\cdot) = E^k[Q(Z, \cdot)]$ satisfies

$$(2.2) \qquad A_k(e_k) = 1/2 - (\sigma/2)\sqrt{(\log M)/n}, \qquad A_k(e_j) = 1/2, \ k \neq j.$$

Therefore $A_k$ admits a unique minimum over the set of vertices $\{e_1, \ldots, e_M\}$ and the minimum is attained at $e_k$.

PROPOSITION 2.1. *Let $Q$ be the squared loss function (2.1). Assume that we observe i.i.d. random vectors $Z_1, \ldots, Z_n$ with the same distribution as $Z$. Denote by $E_n^k$ the expectation with respect to the sample $Z_1, \ldots, Z_n$ when $Z$ has distribution $P_k$. Then there exists an absolute constant $c > 0$ such that*

$$(2.3) \qquad \inf_{T_n} \sup_{k=1,\ldots,M} \left\{ E_n^k[A_k(T_n)] - \min_{1 \leq j \leq M} A_k(e_j) \right\} \geq c\sigma \sqrt{\frac{\log M}{n}},$$

*where the infimum is taken over all the selectors $T_n$.*

A weaker result of similar type [with the rate $1/\sqrt{n}$ instead of $\sqrt{(\log M)/n}$] is given in [14]. Proposition 2.1 implies that the slow rate $\sqrt{(\log M)/n}$ is the best attainable rate for selectors, since the standard ERM selector satisfies the oracle inequality (1.1) with rate $\Delta_{n,M} \sim \sqrt{(\log M)/n}$. Proof of Proposition 2.1 is given in Section 6.

The squared loss function (2.1) satisfies the assumptions of Theorems 4.1 and 4.2 below. As a consequence, the corresponding aggregated estimate $\hat{\theta}_n$, provided by the algorithm of Section 3, attains the bound with fast rate $(\log M)/n$:

$$E_n^k A_k(\hat{\theta}_n) \leq \min_{1 \leq j \leq M} A_k(e_j) + C \frac{(\sigma^2 + 1)\log M}{n} \qquad \forall k = 1, \dots, M.$$

On the other hand, for the same squared loss, Proposition 2.1 shows that a selector with values in $\{e_1, \dots, e_M\}$, in particular the ERM, cannot satisfy an oracle inequality of the type (1.1) with the rate faster than $\sqrt{(\log M)/n}$. This observation suggests that extending the set of possible values of the estimator to the whole simplex $\Theta$ may help to obtain faster rates of aggregation.

**3. The algorithm.** Procedures with values in $\Theta$, that is, convex mixtures of the initial estimators, can be constructed in various ways. One of them originates from the idea of mirror descent due to Nemirovski and Yudin [19]. This idea has been further developed in [3, 20], mainly in the deterministic optimization framework. A version of the mirror descent method due to Nesterov [20] has been applied to the aggregation problem in [12] under the name of mirror averaging. As shown in [12], for convex loss functions $Q$ the mirror averaging estimator $\tilde{\theta}_n$ satisfies under mild assumptions the following oracle inequality:

$$(3.1) \qquad E_n A(\tilde{\theta}_n) \leq \min_{\theta \in \Theta} A(\theta) + C_0 \sqrt{\frac{\log M}{n}},$$

where $C_0 > 0$ is a constant depending only on the supremum norm of the gradient $\nabla_\theta Q(\cdot, \cdot)$. The name mirror averaging reflects the fact that the algorithm does a stochastic gradient descent in the dual space with further "mirroring" to the primal space and averaging; for more details and discussion see [12].

Note that in (3.1) the minimum is taken over the whole simplex $\Theta$, so an inequality of the type (1.1) holds as well, but for large $n$ the remainder term in (3.1) is of larger order than the optimal one given in (1.2).

To improve upon this, consider the following auxiliary stochastic linear programming problem. If $A$ is a convex function, we can bound it from above by a linear function:

$$A(\theta) \leq \sum_{j=1}^{M} \theta^{(j)} A(e_j) \triangleq \tilde{A}(\theta) \qquad \forall \theta \in \Theta,$$

where $\tilde{A}(\theta) = E\tilde{Q}(Z, \theta)$, with

$$\tilde{Q}(Z, \theta) \triangleq \theta^\top u(Z), \qquad u(Z) \triangleq (Q(Z, e_1), \dots, Q(Z, e_M))^\top.$$

Note that

$$\tilde{A}(e_j) = A(e_j), \qquad j = 1, \ldots, M.$$

Since $\Theta$ is a simplex, the minimum of the linear function $\tilde{A}$ is attained at one of its vertices. Therefore,

$$\min_{\theta \in \Theta} \tilde{A}(\theta) = \min_{1 \le j \le M} A(e_j),$$

which shows that the *linear* stochastic programming problem of minimization of $\tilde{A}$ on $\Theta$ is linked to the problem of MS aggregation. This also suggests that the mirror averaging algorithm of [12] applied to minimization of the linear function $\tilde{A}$ could make sense to achieve our MS aggregation goal. Particularizing the definition of mirror averaging procedure from [12] for linear function $\tilde{A}$, we get the following algorithm.

For $\beta > 0$ define the function $W_\beta : \mathbb{R}^M \to \mathbb{R}$ by

$$(3.2) \qquad W_\beta(z) \triangleq \beta \log \left( \frac{1}{M} \sum_{j=1}^M e^{-z^{(j)}/\beta} \right), \qquad z = (z^{(1)}, \ldots, z^{(M)}).$$

The gradient of $W_\beta$ is given by

$$\nabla W_\beta(z) = \left[ -\frac{e^{-z^{(j)}/\beta}}{\sum_{k=1}^M e^{-z^{(k)}/\beta}} \right]_{j=1}^M.$$

Consider the vector

$$u_i \triangleq (Q(Z_i, e_1), \ldots, Q(Z_i, e_M))^\top = u(Z_i) = \nabla_\theta \tilde{Q}(Z_i, \theta),$$

and the iterations:

- Fix the initial values $\theta_0 \in \Theta$ and $\zeta_0 = 0 \in \mathbb{R}^M$.
- For $i = 1, \ldots, n-1$, do the recursive update

$$(3.3) \qquad \begin{aligned} \zeta_i &= \zeta_{i-1} + u_i, \\ \theta_i &= -\nabla W_\beta(\zeta_i). \end{aligned}$$

- Output at iteration $n$ the average

$$(3.4) \qquad \hat{\theta}_n = \frac{1}{n} \sum_{i=1}^n \theta_{i-1}.$$

Note that the estimator $\hat{\theta}_n$ is measurable with respect to $(Z_1, \ldots, Z_{n-1})$. The components $\theta_i^{(j)}$ of the vector $\theta_i$ from (3.3) can be written in the form

$$\theta_i^{(j)} = \frac{\exp\left( -\beta^{-1} \sum_{m=1}^i Q(Z_m, e_j) \right)}{\sum_{k=1}^M \exp\left( -\beta^{-1} \sum_{m=1}^i Q(Z_m, e_k) \right)}, \qquad j = 1, \ldots, M.$$

The "mirroring" function $\nabla W_\beta$ maps the variables $\zeta_i$ that take on values in the dual space (which is $\mathbb{R}^M$ equipped with the $\ell_\infty$ norm) to the primal space (which is the $\ell_1$ body $\Theta$); cf. [12]. Note that $W_\beta$ defined in (3.2) is not the only possible choice; other functions $W_\beta$ satisfying the conditions described in [12] can be used to construct the updates (3.3).

We arrived at the algorithm (3.3)–(3.4) by a linear stochastic programming argument. It is interesting that several particular cases or versions of this algorithm are well known, and they were derived from different considerations. We mention first the literature on prediction of individual deterministic sequences. For a detailed account on this subject see [10]. A general problem considered there is for an agent to compete against the observed predictions of a group of experts, so that the agent's error is close to that of the best expert. In that framework the observations $Z_i$ are supposed to be uniformly bounded nonrandom variables, and the risk function is defined as the cumulative loss over the trajectory. Interestingly, for such problems, which are quite different from ours, methods similar to (3.3) constitute one of the principal tools; cf. [11, 13, 23, 26, 27]. However, in contrast to our procedure, those methods do not involve the averaging step (3.4); they do not need it because they deal with non-random observations and cumulative losses. Note that the algorithm with the averaging step (3.4) included, that is, the one that we consider here, has also been discussed in the literature, though only for two specific combinations of loss function/model: the squared loss $Q$ in regression model [5, 8, 9] and the Kullback–Leibler loss $Q$ in density estimation [7, 9, 29]. It is interesting that in the latter case the algorithm (3.3)–(3.4) can be derived using information-theoretical arguments from the theory of source coding; cf. [9].

Remark that we define algorithm (3.3)–(3.4) for a general loss function $Q$, and we consider arbitrary i.i.d. data $Z_i$, not restricted to a particular model.

Since (3.3)–(3.4) is a special case of the mirror averaging method of [12] corresponding to a linear function $\tilde{A}$, the coarse oracle inequality (3.1) remains valid with $A$ replaced by $\tilde{A}$. But we show below that in fact $\hat{\theta}_n$ satisfies a stronger inequality, that is, one with the optimal rate (1.2).

**4. Main results.** In this section we prove two theorems. They establish oracle inequalities of the type (1.1) for $\hat{\theta}_n$. Theorem 4.1 requires a more conservative assumption on the loss functions $Q$ than Theorem 4.2. This assumption is easier to check, and it often leads to a sharper bound but not for such models as nonparametric density estimation with the $L_2$ loss which will be treated using Theorem 4.2. In some cases (e.g., in regression with Gaussian noise) Theorem 4.1 yields a suboptimal remainder term, while Theorem 4.2 does the correct job. In both theorems it is supposed that the values $A(e_1), \ldots, A(e_M)$ are finite. We will also need the following definition.

DEFINITION 4.1.   A function $T : \mathbb{R}^M \to \mathbb{R}$ is called *exponentially concave* if the composite function $\exp \circ T$ is concave.

It is straightforward to see that exponential concavity of a function $-T$ implies that $T$ is convex. Furthermore, if $-T/\beta$ is exponentially concave for some $\beta > 0$, then $-T/\beta'$ is exponentially concave for all $\beta' > \beta$. Let $Q_1$ be the function on $\mathcal{Z} \times \Theta \times \Theta$ defined by $Q_1(z, \theta, \theta') = Q(z, \theta) - Q(z, \theta')$ for all $z \in \mathcal{Z}$ and all $\theta, \theta' \in \Theta$.

THEOREM 4.1.   *Assume that $Q_1$ can be decomposed into the sum of two functions $Q_1 = Q_2 + Q_3$ such that:*

- *The mapping $\theta \mapsto -Q_2(z, \theta, \theta')/\beta$ is exponentially concave on the simplex $\Theta$, for all $z \in \mathcal{Z}$, $\theta' \in \Theta$, and $Q_2(z, \theta, \theta) = 0$ for all $z \in \mathcal{Z}, \theta \in \Theta$.*
- *There exists a function $R$ on $\mathcal{Z}$ integrable with respect to $P$ and such that $-Q_3(z, \theta, \theta') \leq R(z)$, for all $z \in \mathcal{Z}, \theta, \theta' \in \Theta$.*

*Then the aggregate $\hat{\theta}_n$ satisfies, for any $M \geq 2, n \geq 1$, the following oracle inequality:*

$$E_{n-1} A(\hat{\theta}_n) \leq \min_{1 \leq j \leq M} A(e_j) + \frac{\beta \log M}{n} + E[R(Z)].$$

THEOREM 4.2.   *Assume that for some $\beta > 0$ there exists a Borel function $\Psi_\beta : \Theta \times \Theta \to \mathbb{R}_+$ such that the mapping $\theta \mapsto \Psi_\beta(\theta, \theta')$ is concave on the simplex $\Theta$ for any fixed $\theta' \in \Theta$, $\Psi_\beta(\theta, \theta) = 1$ and $E \exp(-Q_1(Z, \theta, \theta')/\beta) \leq \Psi_\beta(\theta, \theta')$ for all $\theta, \theta' \in \Theta$. Then the aggregate $\hat{\theta}_n$ satisfies, for any $M \geq 2, n \geq 1$, the following oracle inequality:*

$$E_{n-1} A(\hat{\theta}_n) \leq \min_{1 \leq j \leq M} A(e_j) + \frac{\beta \log M}{n}.$$

Proofs of both theorems are based on the following lemma. Introduce the discrete random variable $\omega$ with values in the set $\{e_1, \ldots, e_M\}$ and with the distribution $\mathbb{P}$ defined conditionally on $(Z_1, \ldots, Z_{n-1})$ by $\mathbb{P}[\omega = e_j] = \hat{\theta}_n^{(j)}$ where $\hat{\theta}_n^{(j)}$ is the $j$th component of $\hat{\theta}_n$. The expectation corresponding to $\mathbb{P}$ is denoted by $\mathbb{E}$.

LEMMA 4.1.   *For any measurable function $Q$ and any $\beta > 0$ we have*

$$(4.1) \qquad E_{n-1} A(\hat{\theta}_n) \leq \min_{1 \leq j \leq M} A(e_j) + \frac{\beta \log M}{n} + S_1,$$

*where*

$$S_1 \triangleq \beta E_n \log\left( \mathbb{E} \exp\left[ -\frac{Q_1(Z_n, \omega, \mathbb{E}[\omega])}{\beta} \right] \right).$$

PROOF. By definition of $W_\beta(\cdot)$, for $i = 1, \ldots, n$,

$$
\begin{aligned}
W_\beta(\zeta_i) - W_\beta(\zeta_{i-1}) &= \beta \log\left( \frac{\sum_{j=1}^M e^{-\zeta_i^{(j)}/\beta}}{\sum_{j=1}^M e^{-\zeta_{i-1}^{(j)}/\beta}} \right) \\
&= \beta \log(-v_i^\top \nabla W_\beta(\zeta_{i-1})) = \beta \log(v_i^\top \theta_{i-1}),
\end{aligned}
$$
(4.2)

*where*

$$v_i = \left[ \exp\left( -\frac{u_i^{(j)}}{\beta} \right) \right]_{j=1}^M.$$

Taking expectations on both sides of (4.2), summing up over $i$, using the fact that $(\theta_{i-1}, Z_i)$ has the same distribution as $(\theta_{i-1}, Z_n)$ for $i = 1, \ldots, n$, and applying the Jensen inequality, we get

$$
\begin{aligned}
\frac{E_n[W_\beta(\zeta_n) - W_\beta(\zeta_0)]}{n} &= \frac{\beta}{n} \sum_{i=1}^n E_n \log\left( \sum_{j=1}^M \theta_{i-1}^{(j)} \exp\left[ -\frac{Q(Z_i, e_j)}{\beta} \right] \right) \\
&= \frac{\beta}{n} \sum_{i=1}^n E_n \log\left( \sum_{j=1}^M \theta_{i-1}^{(j)} \exp\left[ -\frac{Q(Z_n, e_j)}{\beta} \right] \right) \\
&\leq \beta E_n \log\left( \sum_{j=1}^M \hat{\theta}_n^{(j)} \exp\left[ -\frac{Q(Z_n, e_j)}{\beta} \right] \right) \triangleq S.
\end{aligned}
$$
(4.3)

Since $Q_1(z, \omega, \mathbb{E}[\omega]) = Q(z, \omega) - Q(z, \mathbb{E}[\omega])$ and $\mathbb{E}[\omega] = \hat{\theta}_n$, the RHS of (4.3) can be written in the form

$$
\begin{aligned}
S &= \beta E_n \log\left( \mathbb{E} \exp\left[ -\frac{Q(Z_n, \omega)}{\beta} \right] \right) \\
&= \beta E_n \log\left( \exp\left[ -\frac{Q(Z_n, \mathbb{E}[\omega])}{\beta} \right] \right) + S_1 \\
&= -E_{n-1} A(\hat{\theta}_n) + S_1.
\end{aligned}
$$
(4.4)

We now bound from below the LHS of (4.3). For any $j^\star = 1, \ldots, M$, by monotonicity of the function $\log(\cdot)$, we have

$$W_\beta(\zeta_n) \geq \beta \log\left( \frac{1}{M} e^{-\zeta_n^{(j^\star)}/\beta} \right) = -\beta \log M - \zeta_n^{(j^\star)},$$

where $\zeta_n^{(j^\star)} = \zeta_n^\top e_{j^\star}$ is the $j^\star$th component of $\zeta_n$. Set $j^\star = \arg\min_{1 \le j \le M} A(e_j)$. Then, using the fact that $W_\beta(\zeta_0) = W_\beta(0) = 0$ we obtain

$$
\begin{aligned}
(4.5) \qquad \frac{E_n[W_\beta(\zeta_n) - W_\beta(\zeta_0)]}{n} &\ge -\frac{\beta \log M}{n} - \frac{E_n[\zeta_n^\top e_{j^\star}]}{n} \\
&= -\frac{\beta \log M}{n} - \min_{1 \le j \le M} A(e_j).
\end{aligned}
$$

Combining (4.3), (4.4) and (4.5) gives the lemma. □

In view of Lemma 4.1, to prove Theorems 4.1 and 4.2 it remains to give appropriate upper bounds for $S_1$.

PROOF OF THEOREM 4.1. Since $Q_1 = Q_2 + Q_3$, with $-Q_3(z, \theta, \theta') \le R(z)$ for all $z \in \mathcal{Z}, \theta, \theta' \in \Theta$, the quantity $S_1$ can be bounded from above as follows:

$$
S_1 \le \beta E_n \log\left(\mathbb{E} \exp\left[-\frac{Q_2(Z_n, \omega, \mathbb{E}[\omega])}{\beta}\right]\right) + E_n[R(Z_n)].
$$

Now since $-Q_2(z, \cdot)/\beta$ is exponentially concave on $\Theta$ for all $z \in \mathcal{Z}$, the Jensen inequality yields

$$
\mathbb{E} \exp\left[-\frac{Q_2(Z_n, \omega, \mathbb{E}[\omega])}{\beta}\right] \le \exp\left[-\frac{Q_2(Z_n, \mathbb{E}[\omega], \mathbb{E}[\omega])}{\beta}\right] = 1.
$$

Therefore $S_1 \le E_n[R(Z_n)]$. This and Lemma 4.1 imply the result of the theorem. □

PROOF OF THEOREM 4.2. Using the Jensen inequality twice, with the concave functions $\log(\cdot)$ and $\Psi_\beta(\cdot, \mathbb{E}[\omega])$, we get

$$
\begin{aligned}
(4.6) \qquad S_1 &\le \beta E_{n-1} \log\left(E\mathbb{E} \exp\left[-\frac{Q_1(Z, \omega, \mathbb{E}[\omega])}{\beta}\right]\right) \\
&= \beta E_{n-1} \log\left(\mathbb{E} E \exp\left[-\frac{Q_1(Z, \omega, \mathbb{E}[\omega])}{\beta}\right]\right) \\
&\le \beta E_{n-1} \log(\mathbb{E}\Psi_\beta(\omega, \mathbb{E}[\omega])) \\
&\le \beta E_{n-1} \log(\Psi_\beta(\mathbb{E}[\omega], \mathbb{E}[\omega])) = 0,
\end{aligned}
$$

where the first equality is due to the Fubini theorem. Theorem 4.2 follows now from (4.6) and Lemma 4.1. □

REMARK. A particular case of Theorem 4.1 where $Q_3 \equiv 0$ and the loss $Q$ is uniformly bounded in $z, \theta$ can be derived from the theory of prediction of deterministic sequences discussed in Section 3 above. We sketch here the argument that can be used. If written in our notation, some results of that theory (see, e.g., [13, 23] or Section 3.3 of [10]) are as follows: under exponential concavity of $\theta \mapsto -\eta Q(z, \theta)$ for some $\eta > 0$ and boundedness of $\sup_{z,\theta} |Q(z, \theta)|$, for any fixed sequence $Z_i$ we have

$$(4.7) \qquad \frac{1}{n} \sum_{i=1}^{n} Q(Z_i, \theta_{i-1}) \leq \min_{j=1,\ldots,M} \frac{1}{n} \sum_{i=1}^{n} Q(Z_i, e_j) + \frac{C \log M}{n}$$

where $C$ is a constant depending only on $\beta$ and on the value $\sup_{z,\theta} |Q(z, \theta)|$. Assuming now that $Z_i$ are random and i.i.d., taking expectations in (4.7) and interchanging the expectation and the minimum on the right-hand side we obtain

$$(4.8) \qquad E_n \left( \frac{1}{n} \sum_{i=1}^{n} Q(Z_i, \theta_{i-1}) \right) \leq \min_{j=1,\ldots,M} A(e_j) + \frac{C \log M}{n}.$$

Now, exponential concavity of $\theta \mapsto -\eta Q(z, \theta)$ implies convexity of $\theta \mapsto Q(z, \theta)$ and thus convexity of $A(\cdot)$. Therefore, since $\theta_{i-1}$ is measurable with respect to $Z_1, \ldots, Z_{i-1}$ using Jensen's inequality and the definition of $\hat{\theta}_n$ we get

$$
\begin{aligned}
(4.9) \quad E_n \left( \frac{1}{n} \sum_{i=1}^{n} Q(Z_i, \theta_{i-1}) \right) &= \frac{1}{n} \sum_{i=1}^{n} E_{i-1} A(\theta_{i-1}) \\
&= E_{n-1} \left( \frac{1}{n} \sum_{i=1}^{n} A(\theta_{i-1}) \right) \geq E_{n-1} A(\hat{\theta}_n).
\end{aligned}
$$

Combining (4.8) and (4.9) we get inequality of the form (1.1)–(1.2). We note that such an argument can be used as an alternate proof of Corollary 5.3 in the next section. However, it does not apply to other examples that we treat below using Theorems 4.1 and 4.2 since in those examples either the loss is not bounded or the exponential concavity condition is not satisfied. We need only some approximate exponential concavity (when using Theorem 4.1) or a kind of "exponential concavity in the mean" (when using Theorem 4.2).

**5. Examples.** In this section we apply Theorems 4.1 and 4.2 to three common statistical problems (regression, classification and density estimation) in order to establish some new oracle inequalities. In particular, we cover the two examples for which our algorithm has been already studied in the literature: regression model with squared loss and density estimation with KL loss. For the latter case we observe that our general argument easily implies the earlier results [7, 9, 29], while for regression with squared loss we significantly improve what was known before [5, 8, 9].

All the loss functions considered below are twice differentiable. The following proposition gives a simple sufficient condition for exponential concavity.

PROPOSITION 5.1. *Let $g$ be a twice differentiable function on $\Theta$ with gradient $\nabla g(\theta)$ and Hessian matrix $\nabla^2 g(\theta)$, $\theta \in \Theta$. If there exists $\beta > 0$ such that for any $\theta \in \Theta$, the matrix*

$$\beta \nabla^2 g(\theta) - \nabla g(\theta)(\nabla g(\theta))^\top,$$

*is positive semidefinite, then $-g(\cdot)/\beta$ is exponentially concave on the simplex $\Theta$.*

PROOF. Since $g$ is twice differentiable $\exp(-g(\cdot)/\beta)$ is also twice differentiable with Hessian matrix

$$(5.1) \qquad \mathcal{H}(\theta) = \frac{1}{\beta} \exp\left(-\frac{g(\theta)}{\beta}\right) \left[ \frac{\nabla g(\theta)(\nabla g(\theta))^\top}{\beta} - \nabla^2 g(\theta) \right].$$

For any $\lambda \in \mathbb{R}^M$, $\theta \in \Theta$, we have

$$\lambda^\top \mathcal{H}(\theta) \lambda = \frac{1}{\beta} \exp\left(-\frac{g(\theta)}{\beta}\right) \left[ \frac{(\lambda^\top \nabla g(\theta))^2}{\beta} - \lambda^\top [\nabla^2 g(\theta)] \lambda \right] \leq 0.$$

Hence $\exp(-g(\cdot)/\beta)$ has a negative semidefinite Hessian and is therefore concave. □

5.1. *Application of Theorem 4.1.* We begin with the models that satisfy assumptions of Theorem 4.1.

1. *Regression with squared loss.* Let $\mathcal{Z} = \mathcal{X} \times \mathbb{R}$ where $\mathcal{X}$ is a complete separable metric space equipped with its Borel $\sigma$-algebra. Consider a random variable $Z = (X, Y)$ with $X \in \mathcal{X}$ and $Y \in \mathbb{R}$. Assume that the conditional expectation $f(X) = E(Y|X)$ exists and define $\xi = Y - E(Y|X)$, so that

$$(5.2) \qquad\qquad Y = f(X) + \xi,$$

where $X \in \mathcal{X}$ is a random variable with probability distribution $P_X$, $Y \in \mathbb{R}$, $f : \mathcal{X} \to \mathbb{R}$ is the regression function and $\xi$ is a real-valued random variable satisfying $E(\xi|X) = 0$. Assume that $E(Y^2) < \infty$ and $\|f\|_\infty \leq L$ for some finite constant $L > 0$ where $\| \cdot \|_\infty$ denotes the $L_\infty(P_X)$-norm. We have $M$ functions $f_1, \ldots, f_M$ such that $\|f_j\|_\infty \leq L, j = 1, \ldots, M$. Define $\|f\|_{2,P_X}^2 = \int_\mathcal{X} f^2(x) P_X(dx)$. Our goal is to construct an aggregate that mimics the oracle $\min_{1 \leq j \leq M} \|f_j - f\|_{2,P_X}^2$. The aggregate is based on the i.i.d. sample $(X_1, Y_1), \ldots, (X_n, Y_n)$ where $(X_i, Y_i)$ have the same distribution as $(X, Y)$. For this model, with $z = (x, y) \in \mathcal{X} \times \mathbb{R}$, define the loss function

$$Q(z, \theta) = (y - \theta^\top H(x))^2 \qquad \forall \theta \in \Theta,$$

with $H(x) = (f_1(x), \ldots, f_M(x))^\top$. It yields for all $z \in \mathcal{Z}, \theta, \theta' \in \Theta$,

$$Q_1(z, \theta, \theta') = Q(z, \theta) - Q(z, \theta') = 2y(\theta' - \theta)^\top H(x) + [\theta^\top H(x)]^2 - [\theta'^\top H(x)]^2.$$

Consider positive constants $b$ and $B$ and assume that $\beta > (b/B)^2$. We now decompose $Q_1$ into the sum $Q_1 = Q_2 + Q_3$, where

$$Q_2(z, \theta, \theta') = 2y\mathbb{1}_{\{|y| < B\beta\}}(\theta' - \theta)^\top H(x) + [\theta^\top H(x)]^2 - [\theta'^\top H(x)]^2$$

$$+ \frac{y^2}{B\beta}[(\theta' - \theta)^\top H(x)]^2 \mathbb{1}_{\{b\sqrt{\beta} < |y| < B\beta\}}$$

and

$$Q_3(z, \theta, \theta') = 2y\mathbb{1}_{\{|y| \geq B\beta\}}(\theta' - \theta)^\top H(x) - \frac{y^2}{B\beta}[(\theta' - \theta)^\top H(x)]^2 \mathbb{1}_{\{b\sqrt{\beta} < |y| < B\beta\}}.$$

We have

$$(5.3) \quad -Q_3(z, \theta, \theta') \leq 4L|y|\mathbb{1}_{\{|y| \geq B\beta\}} + \frac{4L^2 y^2}{B\beta}\mathbb{1}_{\{b\sqrt{\beta} < |y| < B\beta\}} \triangleq R_\beta(y).$$

On the other hand, $Q_2(z, \theta, \theta) = 0, \forall \theta \in \Theta, z \in \mathcal{Z}$ and we can prove that the mapping $\theta \mapsto -Q_2(z, \theta, \theta')/\beta$ is exponentially concave for any $z \in \mathcal{Z}, \theta' \in \Theta$ when $b$ and $B$ are properly chosen. For all $\theta \in \Theta$ and $z = (x, y)$ the gradient and Hessian of $Q_2$ are respectively given by

$$\nabla_\theta Q_2 = \nabla_\theta Q_2(z, \theta, \theta')$$

$$= -2(y\mathbb{1}_{\{|y| < B\beta\}} - \theta^\top H(x))H(x)$$

$$- 2\frac{y^2}{B\beta}\mathbb{1}_{\{b\sqrt{\beta} < |y| < B\beta\}}[(\theta' - \theta)^\top H(x)]H(x)$$

and

$$\nabla^2_{\theta\theta} Q_2 = \nabla^2_{\theta\theta} Q_2(z, \theta, \theta') = 2H(x)H(x)^\top + 2\frac{y^2}{B\beta}\mathbb{1}_{\{b\sqrt{\beta} < |y| < B\beta\}}H(x)H(x)^\top.$$

We now prove that Proposition 5.1 applies for $g(\theta) = Q_2(z, \theta, \theta')$, for all $z = (x, y) \in \mathcal{Z}$ and $\theta' \in \Theta$. For any $\lambda \in \mathbb{R}^M$, any $\theta, \theta' \in \Theta$ and any $z \in \mathcal{Z}$,

$$(\lambda^\top \nabla_\theta Q_2)^2 \leq \left(2|y|\mathbb{1}_{\{|y| < B\beta\}} + 2L + \frac{4Ly^2}{B\beta}\mathbb{1}_{\{b\sqrt{\beta} < |y| < B\beta\}}\right)^2 [\lambda^\top H(x)]^2.$$

Note now that $|y| \leq B\beta$ implies that $y^2/B\beta \leq |y|$. Hence

$$(\lambda^\top \nabla_\theta Q_2)^2 \leq (2|y|\mathbb{1}_{\{|y| \leq b\sqrt{\beta}\}} + 2L + (4L + 2)|y|\mathbb{1}_{\{b\sqrt{\beta} < |y| < B\beta\}})^2 [\lambda^\top H(x)]^2$$

$$\leq (8b^2\beta + 8L^2 + 2(4L + 2)^2|y|^2\mathbb{1}_{\{b\sqrt{\beta} < |y| < B\beta\}})[\lambda^\top H(x)]^2.$$

Therefore

$$\frac{(\lambda^\top \nabla_\theta Q_2)^2}{\beta} - \lambda^\top(\nabla^2_{\theta\theta} Q_2)\lambda$$

$$\leq \left(8b^2 + \frac{8L^2}{\beta} - 2 + \left[2(4L + 2)^2 - \frac{2}{B}\right]\frac{|y|^2}{\beta}\mathbb{1}_{\{b\sqrt{\beta} < |y| < B\beta\}}\right)[\lambda^\top H(x)]^2.$$

If we choose $B \leq (4L+2)^{-2}$ and $LB < b < 1/4$, the above quadratic form is smaller than or equal to 0 and Proposition 5.1 applies for any $\beta > (b/B)^2$. Now, since $A(\theta) = EQ(Z, \theta) = E(Y - \theta^\top H(X))^2 = \|f - \theta^\top H\|_{2,P_X}^2 + E(\xi^2)$ for all $\theta \in \Theta$, we obtain the following corollary of Theorem 4.1.

COROLLARY 5.1. *Consider the regression model* (5.2) *where* $X \in \mathcal{X}$, $Y \in \mathbb{R}$, $f : \mathcal{X} \to \mathbb{R}$ *and* $\xi = Y - f(X)$ *is a real-valued random variable satisfying* $E(\xi|X) = 0$. *Assume also that* $E(Y^2) < \infty$ *and* $\|f_j\|_\infty \leq L, j = 1, \ldots, M$, *for some finite constant* $L > 0$. *Then for any positive constants* $B \geq (4L + 2)^{-2}, LB < b < 1/4$ *and any* $\beta \geq (b/B)^2$, *the aggregate estimator* $\tilde{f}_n(x) = \hat{\theta}_n^\top H(x), x \in \mathcal{X}$, *where* $\hat{\theta}_n$ *is obtained by the mirror averaging algorithm, satisfies*

$$(5.4) \quad E_{n-1}\|\tilde{f}_n - f\|_{2,P_X}^2 \leq \min_{1 \leq j \leq M} \|f_j - f\|_{2,P_X}^2 + \frac{\beta \log M}{n} + E[R_\beta(Y)],$$

*where*

$$R_\beta(y) = 4L|y| \mathbb{1}_{\{|y| \geq B\beta\}} + \frac{4L^2 y^2}{B\beta} \mathbb{1}_{\{b\sqrt{\beta} < |y| < B\beta\}}.$$

This result improves an inequality obtained by [5]: it yields better rate under the same moment conditions. We note that the aggregate $\tilde{f}_n$ as in Corollary 5.1 is of the form suggested by Catoni [8, 9]. If there exists a constant $L_0 > 0$ such that $|Y| \leq L_0$ a.s., the last summand disappears for $\beta > 16L_0^2$, and in this case (5.4) can be also deduced from [8, 9], though in a coarser form and under a more restrictive assumption on $\beta$.

An advantage of Corollary 5.1 is that no heavy assumption on the moments of $\xi$ is needed to get reasonable bounds. Thus, the second moment assumption on $Y$ is enough for a bound with the $n^{-1/2}$ rate. Indeed, choosing $\beta \sim (n/\log M)^{2/(2+s)}$, $s > 0$, in Corollary 5.1, we immediately get the following result.

COROLLARY 5.2. *Consider the regression model* (5.2) *where* $X \in \mathcal{X}$, $Y \in \mathbb{R}$, $f : \mathcal{X} \to \mathbb{R}$ *and* $\xi = Y - f(X)$ *is a real-valued random variable satisfying* $E(\xi|X) = 0$. *Assume also that* $E(|Y|^s) \leq m_s < \infty$ *for some* $s \geq 2$ *and* $\|f_j\|_\infty \leq L, j = 1, \ldots, M$, *for some finite constant* $L > 0$. *Then there exist constants* $C_1 > 0$ *and* $C_2 = C_2(m_s, L, C_1) > 0$ *such that the aggregate estimator* $\tilde{f}_n(x) = \hat{\theta}_n^\top H(x), x \in \mathcal{X}$, *where* $\hat{\theta}_n$ *is obtained by the mirror averaging algorithm with* $\beta = C_1(n/\log M)^{2/(2+s)}$, *satisfies*

$$(5.5) \quad E_{n-1}\|\tilde{f}_n - f\|_{2,P_X}^2 \leq \min_{1 \leq j \leq M} \|f_j - f\|_{2,P_X}^2 + C_2 \left(\frac{\log M}{n}\right)^{s/(2+s)}.$$

2. *Classification.* Consider the problem of binary classification. Let $(\mathcal{X}, \mathcal{F})$ be a measurable space, and set $\mathcal{Z} = \mathcal{X} \times \{-1, 1\}$. Consider $Z = (X, Y)$ where $X$ is a random variable with values in $\mathcal{X}$ and $Y$ is a random label with values in $\{-1, 1\}$. For a fixed convex twice differentiable function $\varphi \colon \mathbb{R} \to \mathbb{R}_+$, define the $\varphi$-risk of a real-valued classifier $h \colon \mathcal{X} \to [-1, 1]$ as $E\varphi(-Yh(X))$. In our framework, we have $M$ such classifiers $h_1, \ldots, h_M$ and the goal is to mimic the oracle $\min_{1 \le j \le M} E\varphi(-Yh_j(X))$ based on the i.i.d. sample $(X_1, Y_1), \ldots, (X_n, Y_n)$ where $(X_i, Y_i)$ have the same distribution as $(X, Y)$. For any $z = (x, y) \in \mathcal{X} \times \{-1, 1\}$, we define the loss function

$$Q(z, \theta) = \varphi(-y\theta^\top H(x)) \ge 0 \qquad \forall \theta \in \Theta,$$

where $H(x) = (h_1(x), \ldots, h_M(x))^\top$. For such a function and for all $\theta \in \Theta$, $z = (x, y) \in \mathcal{X} \times \{-1, 1\}$ we have

$$\nabla_\theta Q_1(z, \theta, \theta') = -y\varphi'(-y\theta^\top H(x))H(x),$$

$$\nabla^2_{\theta\theta} Q_1(z, \theta, \theta') = \varphi''(-y\theta^\top H(x))H(x)H(x)^\top.$$

Thus, from Proposition 5.1 the mapping $\theta \mapsto -Q_1(z, \theta, \theta')/\beta$ is exponentially concave for all $z$ and $\theta'$ if $\beta \ge \beta_\varphi$ where $\beta_\varphi$ is such that $[\varphi'(x)]^2 \le \beta_\varphi \varphi''(x)$, $\forall |x| \le 1$. Now, since

$$A(\theta) = EQ(Z, \theta) \quad \text{and} \quad Q(Z, \theta) = \varphi(-Y\theta^\top H(X)), \qquad \forall \theta \in \Theta, Z = (X, Y),$$

we obtain the following corollary of Theorem 4.1 applied with $Q_2 = Q_1$ and $Q_3 \equiv 0$.

COROLLARY 5.3. *Consider the binary classification problem as described above. Assume that the convex function $\varphi$ is such that*

$$[\varphi'(x)]^2 \le \beta_\varphi \varphi''(x) \qquad \forall |x| \le 1.$$

*Then the aggregate classifier $\tilde{h}_n(x) = \hat{\theta}_n^\top H(x), x \in \mathcal{X}$, where $\hat{\theta}_n$ is obtained by the mirror averaging algorithm with $\beta \ge \beta_\varphi$, satisfies*

$$(5.6) \qquad E_n \varphi(-Y_n \tilde{h}_n(X_n)) \le \min_{1 \le j \le M} E\varphi(-Y\tilde{h}_j(X)) + \frac{\beta \log M}{n}.$$

For example, inequality (5.6) holds with the exponential Boosting loss $\varphi_1(x) = e^x$, for which $\beta_{\varphi_1} = e$ and for the Logit-Boosting loss $\varphi_2(x) = \log_2(1 + e^x)$ (in that case $\beta_{\varphi_2} = e\log 2$). For the squared loss $\varphi_3(x) = (1 - x)^2$ and the 2-norm soft margin loss $\varphi_4(x) = \max\{0, 1 - x\}^2$ inequality (5.6) is satisfied with $\beta \ge 2$.

3. *Nonparametric density estimation with Kullback–Leibler (KL) loss.* Let $X$ be a random variable with values in a measurable space $(\mathcal{X}, \mathcal{F})$. Assume that the distribution of $X$ admits a density $p$ with respect to a $\sigma$-finite measure $\mu$ on $(\mathcal{X}, \mathcal{F})$. Assume also that we have $M$ probability densities $p_j$ with respect to $\mu$ on $(\mathcal{X}, \mathcal{F})$ (estimators of $p$) and of an i.i.d. sample $X_1, \ldots, X_n$ where $X_i$ take values in $\mathcal{X}$, and have the same distribution as $X$. Define the KL divergence between two probability densities $p$ and $q$ with respect to $\mu$ as

$$\mathcal{K}(p, q) \triangleq \int_{\mathcal{X}} \log\left(\frac{p(x)}{q(x)}\right) p(x) \mu(dx),$$

if the probability distribution corresponding to $p$ is absolutely continuous with respect to the one corresponding to $q$, and $\mathcal{K}(p, q) = \infty$ otherwise. We assume that the entropy integral $\int p(x) \log p(x) \mu(dx)$ is finite.

Our goal is to construct an aggregate that mimics the KL oracle defined by $\min_{1 \le j \le M} \mathcal{K}(p, p_j)$. For $x \in \mathcal{X}$, $\theta \in \Theta$, we introduce the corresponding loss function

(5.7) $$Q(x, \theta) = -\log(\theta^\top H(x)),$$

where $H(x) = (p_1(x), \ldots, p_M(x))^\top$. We set $Z = X$. Then

$$A(\theta) = EQ(X, \theta) = -\int \log(\theta^\top H(x)) p(x) \mu(dx)$$

where the integral is finite if all the divergences $\mathcal{K}(p, p_j)$ are finite. In particular, $A(e_j) = \mathcal{K}(p, p_j) - \int p(x) \log p(x) \mu(dx)$. Since, for all $x \in \mathcal{X}$, we have

$$\exp(-Q_1(x, \theta, \theta')/\beta) = (\theta^\top H(x))^{1/\beta} (\theta'^\top H(x))^{-1/\beta},$$

the mapping $\theta \mapsto -Q_1(x, \theta, \theta')/\beta$ is exponentially concave on $\Theta$ for any $\beta \ge 1$. Hence, we can apply Theorem 4.1, again with $Q_2 = Q_1$ and $Q_3 \equiv 0$ and we obtain the following corollary.

COROLLARY 5.4. *Consider the density estimation problem with the KL loss as described above, such that $\int p(x) |\log p(x)| \mu(dx) < \infty$. Then the aggregate estimator $\tilde{p}_n(x) = \hat{\theta}_n^\top H(x), x \in \mathcal{X}$, where $\hat{\theta}_n$ is obtained by the mirror averaging algorithm with $\beta = 1$, satisfies*

$$E_{n-1} \mathcal{K}(p, \tilde{p}_n) \le \min_{1 \le j \le M} \mathcal{K}(p, p_j) + \frac{\log M}{n}.$$

We note that the KL aggregate $\tilde{p}_n$ as in Corollary 5.4 coincides with the "progressive mixture rule" considered by Catoni [7, 8, 9] and Yang [29] and the oracle inequality of Corollary 5.4 is the one obtained in those papers. We also note that this is the most trivial example of application of our results. In fact, when $Q$ is of the particular form (5.7), the convexity argument that we developed in Theorems 4.1 and 4.2 is not needed since $S_1 = 0$, so that Corollary 5.4 follows directly from Lemma 4.1. Writing the proof of Lemma 4.1 for this particular $Q$ we essentially recover the proof of Theorem 3.1.1 in [9]. Extension of Corollary 5.4 to $\beta \geq 1$ is straightforward but the oracle inequality for the corresponding aggregate ("Gibbs estimator"; cf. [9]) is less interesting because it has obviously a larger remainder term.

5.2. *Applications of Theorem 4.2.* We now apply Theorem 4.2 to obtain bounds for the regression setup that are sharper than the existing ones. We also use this result to handle the problems of density estimation with squared loss and some examples of parametric estimation that cannot be treated using Theorem 4.1.

4. *Regression with squared loss and finite exponential moment.* We consider here the regression model described in Corollary 5.1 under the additional assumption that, conditionally on $X$, the regression residual $\xi$ admits an exponential moment, that is, there exist positive constants $b$ and $D$ such that, $P_X$-a.s.,

$$E(\exp(b|\xi|)|X) \leq D.$$

Since $E(\xi|X) = 0$, this assumption is equivalent to the existence of positive constants $b_0$ and $\sigma^2$ such that, $P_X$-a.s.,

$$(5.8) \qquad E(\exp(t\xi)|X) \leq \exp(\sigma^2 t^2/2) \qquad \forall |t| \leq b_0;$$

cf. [21], page 56.

In this case, application of Corollary 5.1 leads to suboptimal rates because of the term $E[R_\beta(Y)]$ in (5.4). We show now that, using Theorem 4.2, we can obtain an oracle inequality with optimal rate $(\log M)/n$.

To apply Theorem 4.2, we analyze the mapping $\theta \mapsto E \exp(-Q_1(Z, \theta, \theta')/\beta)$. For the regression model with squared loss as described above, we have $Z = (X, Y)$, $Q(Z, \theta) = (Y - \theta^\top H(X))^2$, and

$$E \exp(-Q_1(Z, \theta, \theta')/\beta)$$
$$= E \exp\left(-\frac{1}{\beta}[(Y - H(X)^\top \theta)^2 - (Y - H(X)^\top \theta')^2]\right)$$
$$= E \exp\left(-\frac{1}{\beta}[-2\xi(U(X, \theta) - U(X, \theta')) + U^2(X, \theta) - U^2(X, \theta')]\right),$$

where $U(X, \theta) \triangleq f(X) - H(X)^\top \theta$. Since

$$|2(U(X, \theta) - U(X, \theta'))| = 2|(\theta - \theta')^\top H(X)| \leq 4L,$$

conditioning on $X$ and using (5.8) we get that, for any $\beta \geq 4L/b_0$,

$$E \exp(-Q_1(Z, \theta, \theta')/\beta) \leq \Psi_\beta(\theta, \theta'),$$

where

$$\Psi_\beta(\theta, \theta') \triangleq E \exp\left(\frac{2\sigma^2}{\beta^2}[(\theta - \theta')^\top H(X)]^2 - \frac{1}{\beta}[U^2(X, \theta) - U^2(X, \theta')]\right).$$

Clearly, $\Psi_\beta(\theta, \theta) = 1$. Thus, to apply Theorem 4.2 it suffices now to specify $\beta_0 > 0$ such that the mapping

$$\theta \mapsto \bar{Q}(x, \theta, \theta') \triangleq \left(-\frac{1}{\beta} + \frac{2\sigma^2}{\beta^2}\right)(\theta^\top H(x))^2$$
$$- \frac{4\sigma^2}{\beta^2}(H(x)^\top \theta)(H(x)^\top \theta') + \frac{2}{\beta}f(x)(H(x)^\top \theta)$$

is exponentially concave for all $\beta \geq \beta_0$, $\theta' \in \Theta$ and almost all $x \in \mathcal{X}$. Note that

$$\nabla_\theta \bar{Q}(x, \theta, \theta') = \left(2\gamma(f(x) - H(x)^\top \theta) + \frac{4\sigma^2}{\beta^2}(f(x) - H(x)^\top \theta')\right)H(x),$$
$$\nabla^2_{\theta\theta} \bar{Q}(x, \theta, \theta') = -2\gamma H(x)H(x)^\top,$$

where $\gamma = \frac{1}{\beta} - \frac{2\sigma^2}{\beta^2}$. Proposition 5.1 implies that $\bar{Q}$ is exponentially concave in $\theta$ if $\nabla^2_{\theta\theta}\bar{Q}(x, \theta, \theta') + \nabla_\theta \bar{Q}(x, \theta, \theta')(\nabla_\theta \bar{Q}(x, \theta, \theta'))^\top \leq 0$. If we assume that $\max_j \|f - f_j\|_\infty \leq \tilde{L}$, we obtain that the latter property holds for $\beta \geq \beta_0 \triangleq 2\sigma^2 + 2\tilde{L}^2$. Thus, Theorem 4.2 applies for $\beta \geq \max(2\sigma^2 + 2\tilde{L}^2, 4L/b_0)$ and we have proved the following result.

COROLLARY 5.5. *Consider the regression model (5.2) where $X \in \mathcal{X}$, $Y \in \mathbb{R}$, $f: \mathcal{X} \to \mathbb{R}$ and the random variable $\xi = Y - f(X)$ is such that there exist positive constants $b_0$ and $\sigma^2$ for which (5.8) holds $P_X$-a.s. Assume also that $\|f - f_j\|_\infty \leq \tilde{L}$ and $\|f_j\|_\infty \leq L, j = 1, \ldots, M$, for some finite positive constants $L, \tilde{L}$. Then for any $\beta \geq \max(2\sigma^2 + 2\tilde{L}^2, 4L/b_0)$ the aggregate estimator $\tilde{f}_n(x) = \hat{\theta}_n^\top H(x), x \in \mathcal{X}$, where $\hat{\theta}_n$ is obtained by the mirror averaging algorithm, satisfies*

$$(5.9) \qquad E_{n-1}\|\tilde{f}_n - f\|^2_{2, P_X} \leq \min_{1 \leq j \leq M}\|f_j - f\|^2_{2, P_X} + \frac{\beta \log M}{n}.$$

To see how good the constants are, we may compare this corollary with the results obtained in other papers for the particular case where $\xi$ is conditionally Gaussian given $X$. In this case we have $b_0 = \infty$ and Corollary 5.5 yields the following result.

COROLLARY 5.6. *Consider the regression model (5.2) where $X \in \mathcal{X}$, $Y \in \mathbb{R}$, $f: \mathcal{X} \to \mathbb{R}$ and, conditionally on $X$, the random variable $\xi = Y - f(X)$ is Gaussian with zero mean and variance bounded by $\sigma^2$. Assume that $\|f - f_j\|_\infty \leq \tilde{L}$, for some finite constant $\tilde{L} > 0$. Then for any $\beta \geq 2\sigma^2 + 2\tilde{L}^2$ the aggregate estimator $\tilde{f}_n(x) = \hat{\theta}_n^\top H(x), x \in \mathcal{X}$, where $\hat{\theta}_n$ is obtained by the mirror averaging algorithm, satisfies (5.9).*

This result for Gaussian regression model is more general than that of [9], page 89, because we do not assume that $f$ and all $f_j, j = 1, \ldots, M$, are uniformly bounded. Even if we assume in addition that $f$ and all $f_j, j = 1, \ldots, M$, are uniformly bounded by $L$, Corollary 5.6 improves the result of [9], page 89. Indeed, in this case we have $\tilde{L} \leq 2L$ and a sufficient condition on $\beta$ in Corollary 5.6 is $\beta \geq 2\sigma^2 + 8L^2$. In [9], page 89, we find the result of Corollary 5.6, though under much more restrictive condition $\beta \geq 18.01\sigma^2 + 70.4L^2$.

5. *Nonparametric density estimation with the $L_2$ loss.* Let $\mu$ be a $\sigma$-finite measure on the measurable space $(\mathcal{X}, \mathcal{F})$. In this whole example, densities are understood with respect to $\mu$ and $\|\cdot\|_\infty$ denotes the $L_\infty(\mu)$-norm. Assume that we have $M$ probability densities $p_j, \|p_j\|_\infty \leq L, j = 1, \ldots, M$, and of an i.i.d. sample $X_1, \ldots, X_n$ where $X_i$ take values in $\mathcal{X}$, and are distributed as a random variable $X$ with unknown probability density $p$ such that $\|p\|_\infty \leq L$ for some positive constant $L$. Our goal is to mimic the oracle defined by $\min_{1 \leq j \leq M} \|p_j - p\|_2^2$, where $\|p\|_2^2 = \int p^2(x)\mu(dx)$.

The corresponding loss function is defined, for any $x \in \mathcal{X}, \theta \in \Theta$, by

$$Q(x, \theta) = \theta^\top G \theta - 2\theta^\top H(x), \tag{5.10}$$

where $H(x) = (p_1(x), \ldots, p_M(x))^\top$ and $G$ is an $M \times M$ positive semidefinite matrix with elements $G_{jk} = \int p_j p_k \, d\mu \leq L$. We set $Z = X$. Then $A(\theta) = EQ(X, \theta) = \|p - \theta^\top H\|_2^2 - \|p\|_2^2$. We now want to check conditions of Theorem 4.2, that is, to show that for the loss function (5.10), the mapping $\theta \mapsto E\exp(-Q_1(X, \theta, \theta')/\beta)$ is concave on $\Theta$, for any $\theta' \in \Theta$ and for $\beta \geq \beta_0$ with some $\beta_0 > 0$ that will be specified below. Note first that

$$\begin{aligned}
Q_1(x, \theta, \theta') &= Q(x, \theta) - Q(x, \theta') \\
&= (\theta - \theta')^\top G(\theta + \theta') - 2(\theta - \theta')^\top H(x).
\end{aligned} \tag{5.11}$$

Fix $\theta' \in \Theta$. Concavity of the above mapping can be checked by considering its Hessian $\tilde{\mathcal{H}}$ which, in view of (5.1), satisfies for any $\lambda \in \mathbb{R}^M$, $\theta \in \Theta$,

$$\lambda^\top \tilde{\mathcal{H}}(\theta)\lambda = \frac{1}{\beta^2} E\left\{ \exp\left(-\frac{Q_1(X, \theta, \theta')}{\beta}\right) [(\lambda^\top \nabla_\theta Q_1(X, \theta, \theta'))^2 \right. \\
\left. - \beta \lambda^\top \nabla_{\theta\theta}^2 Q_1(X, \theta, \theta')\lambda] \right\}.$$

Note that for any $x \in \mathcal{X}, \theta \in \Theta$ we have

$$\nabla_\theta Q_1(x, \theta, \theta') = 2G\theta - 2H(x) \quad \text{and} \quad \nabla_{\theta\theta}^2 Q_1(x, \theta, \theta') = 2G.$$

By (5.11) this yields, for any $\lambda \in \mathbb{R}^M$, $\theta, \theta' \in \Theta$,

$$\lambda^\top \tilde{\mathcal{H}}(\theta)\lambda = -\frac{2}{\beta^2} E\left\{ \exp\left(-\frac{(\theta - \theta')^\top G(\theta + \theta') - 2(\theta - \theta')^\top H(X)}{\beta}\right) \right.$$

$$(5.12) \qquad\qquad \times [\beta\lambda^\top G\lambda - 2(\lambda^\top(G\theta - H(X)))^2]\Big\}$$

$$\leq -\frac{2}{\beta^2}\exp\Big(-\frac{(\theta - \theta')^\top G(\theta + \theta')}{\beta}\Big)F(\lambda, \theta, \theta'),$$

where

$$(5.13) \qquad \begin{aligned} F(\lambda, \theta, \theta') = E\Big\{&\exp\Big(\frac{2(\theta - \theta')^\top H(X)}{\beta}\Big) \\ &\times [\beta\lambda^\top G\lambda - 4(\lambda^\top G\theta)^2 - 4(\lambda^\top H(X))^2]\Big\}. \end{aligned}$$

Observe that by the Cauchy inequality

$$(5.14) \qquad (\lambda^\top G\theta)^2 \leq \lambda^\top G\lambda\theta^\top G\theta \leq L\lambda^\top G\lambda \qquad \forall\theta \in \Theta.$$

Further,

$$(5.15) \qquad \begin{aligned} E(\lambda^\top H(X))^2 &= \int(\lambda^\top H(x))^2 p(x)\mu(dx) \\ &\leq L\int(\lambda^\top H(x))^2\mu(dx) = L\lambda^\top G\lambda. \end{aligned}$$

Using (5.14) and (5.15) and the fact that $\|\theta - \theta'\|_1 \leq 2$ where $\|\cdot\|_1$ stands for the $\ell_1(\mathbb{R}^M)$-norm, we obtain

$$\begin{aligned} F(\lambda, \theta, \theta') &\geq (\beta - 4L)\lambda^\top G\lambda E\exp\Big(\frac{2(\theta - \theta')^\top H(X)}{\beta}\Big) \\ &\quad - 4E\Big\{\exp\Big(\frac{2(\theta - \theta')^\top H(X)}{\beta}\Big)(\lambda^\top H(X))^2\Big\} \\ &\geq (\beta - 4L)\lambda^\top G\lambda\exp\Big(-\frac{4L}{\beta}\Big) - 4L\lambda^\top G\lambda\exp\Big(\frac{4L}{\beta}\Big) \geq 0 \end{aligned}$$

provided that

$$\frac{\beta - 4L}{4L}\exp\Big(-\frac{8L}{\beta}\Big) \geq 1.$$

Note that the last inequality is guaranteed for $\beta \geq \beta_0 = 12L$. We conclude that for $\beta \geq 12L$ the Hessian $\tilde{\mathcal{H}}$ in (5.12) is negative semidefinite and therefore the mapping $\theta \mapsto E\exp(-Q_1(X, \theta, \theta')/\beta)$ is concave on $\Theta$ for any fixed $\theta' \in \Theta$. Thus we have proved the following corollary of Theorem 4.2.

COROLLARY 5.7. *Consider the density estimation problem with the $L_2$ loss as described above. Then, for any $\beta \geq 12L$, the aggregate estimator* $\tilde{p}_n(x) = \hat{\theta}_n^\top H(x), x \in \mathcal{X}$, where $\hat{\theta}_n$ is obtained by the mirror averaging algorithm, satisfies

$$E_{n-1}\|\tilde{p}_n - p\|_2^2 \le \min_{1 \le j \le M} \|p_j - p\|_2^2 + \frac{\beta \log M}{n}.$$

6. *Parametric estimation with Kullback–Leibler (KL) loss.* Let $\mathcal{P} = \{P_a, a \in \mathcal{A}\}$ be a family of probability measures on a measurable space $(\mathcal{X}, \mathcal{F})$ dominated by a $\sigma$-finite measure $\mu$ on $(\mathcal{X}, \mathcal{F})$. Here $\mathcal{A} \subset \mathbb{R}^m$ is a bounded set of parameters. The densities relative to $\mu$ are denoted by $p(x, a) = (dP_a/d\mu)(x), x \in \mathcal{X}$. Let $X$ be a random variable with values in $\mathcal{X}$ distributed according to $P_{a^*}$ where $a^* \in \mathcal{A}$ is the unknown true value of the parameter.

In the aggregation framework, we have $M$ values $a_1, \ldots, a_M \in \mathcal{A}$ (preliminary estimators of $a$) and of an i.i.d. sample $X_1, \ldots, X_n$ where $X_i$ take values in $\mathcal{X}$, and have the same distribution as $X$. Our goal is to construct an aggregate $\tilde{a}_n$ that mimics the parametric KL oracle defined by $\min_{1 \le j \le M} K(a^*, a_j)$, where

$$K(a, b) \triangleq \mathcal{K}(p(\cdot, a), p(\cdot, b)) \qquad \forall a, b \in \mathcal{A}.$$

For $x \in \mathcal{X}$, $\theta \in \Theta$, we introduce the corresponding loss function

$$Q(x, \theta) = -\log p(x, \theta^\top H),$$

where $H = (a_1, \ldots, a_M)^\top$. We set $Z = X$. Then

$$A(\theta) = EQ(X, \theta) = -\int \log(p(x, \theta^\top H)) p(x, a^*) \mu(dx),$$

$$A(e_j) = K(a^*, a_j) - \int p(x, a^*) \log(p(x, a^*)) \mu(dx).$$

Since, for all $x \in \mathcal{X}$, $\exp(-Q(x, \theta)/\beta) = (p(x, \theta^\top H))^{1/\beta}$, to apply Theorem 4.2 we need the following assumption.

ASSUMPTION 5.1. For some $\beta > 0$ and for any $a \in \mathcal{A}$ there exists a Borel function $\Psi_\beta : \Theta \times \Theta \to \mathbb{R}_+$ such that $\theta \mapsto \Psi_\beta(\theta, \theta')$ is concave on the simplex $\Theta$ for all $\theta' \in \Theta$, $\Psi_\beta(\theta, \theta) = 1$ and

$$\int \left( \frac{p(x, H^\top \theta)}{p(x, H^\top \theta')} \right)^{1/\beta} p(x, a) \mu(dx) \le \Psi_\beta(\theta, \theta')$$

for all $\theta, \theta' \in \Theta$.

COROLLARY 5.8. *Consider the parametric estimation problem with the KL loss as described above and let $\int p(x, a^*) |\log p(x, a^*)| \mu(dx) < \infty$. Suppose that Assumption 5.1 is fulfilled for some $\beta > 0$. Then the aggregate estimator*

$\tilde{a}_n = \hat{\theta}_n^\top H$ of the parameter $a^*$, where $\hat{\theta}_n$ is obtained by the mirror averaging algorithm, satisfies

$$(5.16) \qquad E_{n-1} K(a^*, \tilde{a}_n) \le \min_{1 \le j \le M} K(a^*, a_j) + \frac{\beta \log M}{n}.$$

Aggregation procedures can be used to construct pointwise adaptive locally parametric estimators in nonparametric regression [2]. In this case inequality (5.16) can be applied to prove the corresponding adaptive risk bounds. We now check that Assumption 5.1 is satisfied for several standard parametric families.

- **Univariate Gaussian distribution.** Let $\mu$ be the Lebesgue measure on $\mathbb{R}$ and let $p(x, a) = (\sigma\sqrt{2\pi})^{-1} \exp(-(x-a)^2/(2\sigma^2))$ be the univariate Gaussian density with mean $a \in \mathcal{A} = [-L, L]$ and known variance $\sigma^2 > 0$. Replacing $f(x)$ by $a^*$ and $H(x)$ by $H$ in the proof of Corollary 5.6, and following exactly the same argument as there we find that Assumption 5.1 is satisfied for any $\beta \ge \beta_0 = 2\sigma^2 + 8L^2$. Hence, (5.16) also holds for such $\beta$. Note that in this case $K(a^*, a) = (a^* - a)^2/(2\sigma^2)$.

- **Bernoulli distribution.** Let $\mu$ be the discrete measure on $\{0, 1\}$ such that $\mu(0) = \mu(1) = 1$ and let $p(x, a) = a\mathbb{1}_{\{x=0\}} + (1-a)\mathbb{1}_{\{x=1\}}$ be the density of a Bernoulli random variable with parameter $a \in \mathcal{A} = (0, 1)$. Then

$$\int \left( \frac{p(x, H^\top\theta)}{p(x, H^\top\theta')} \right)^{1/\beta} p(x, a)\mu(dx)$$
$$= \left( \frac{H^\top\theta}{H^\top\theta'} \right)^{1/\beta} a + \left( \frac{1 - H^\top\theta}{1 - H^\top\theta'} \right)^{1/\beta} (1 - a) \triangleq \Psi_\beta(\theta, \theta').$$

This function is concave in $\theta$ for any $\theta' \in \Theta$ if $\beta \ge 1$ and obviously $\Psi_\beta(\theta, \theta) = 1$. Therefore Assumption 5.1 is satisfied and Corollary 5.8 applies with $\beta = 1$.

- **Poisson distribution.** Let $\mu$ be the counting measure on the set of the nonnegative integers $\mathbb{N}$: $\mu(k) = 1, \forall k \in \mathbb{N}$, and let $p(x, a) = \sum_{k=0}^\infty \frac{a^k}{k!} e^{-a} \mathbb{1}_{\{x=k\}}$ be the density of a Poisson random variable with parameter $a \in \mathcal{A} = [\ell, L]$ where $0 < \ell < L < \infty$. Then

$$(5.17)
\begin{aligned}
&\int \left( \frac{p(x, H^\top\theta)}{p(x, H^\top\theta')} \right)^{1/\beta} p(x, a)\mu(dx) \\
&= \exp\left[ a\left( \frac{H^\top\theta}{H^\top\theta'} \right)^{1/\beta} - a - \frac{H^\top(\theta - \theta')}{\beta} \right] \triangleq \Psi_\beta(\theta, \theta').
\end{aligned}$$

Clearly, $\Psi_\beta(\theta, \theta) = 1$ and it is not hard to show that $\Psi_\beta$ in (5.17) is concave as a function of $\theta$ for any $\theta' \in \Theta$, provided that $\beta \ge 1 + L(1 + L/\ell)(L/\ell)^{1/(2L+1)}$. Therefore Assumption 5.1 is satisfied and Corollary 5.8 applies with $\beta \ge \beta_0 = 1 + L(1 + L/\ell)(L/\ell)^{1/(2L+1)}$.

**6. Proof of Proposition 2.1.** In view of (2.2), for any selector $T_n$ constructed from the observations $Z_1, \ldots, Z_n$ we have

$$A_k(T_n) - \min_{1 \leq j \leq M} A_k(e_j) \geq \left[ A_k(T_n) - \min_{1 \leq j \leq M} A_k(e_j) \right] \mathbb{1}_{\{T_n \neq e_k\}}$$

$$= \sigma \sqrt{\frac{\log M}{n}} \mathbb{1}_{\{T_n \neq e_k\}}.$$

Taking expectation on both sides of the previous inequality yields

$$E_n^k[A_k(T_n)] - \min_{1 \leq j \leq M} A_k(e_j) \geq \sigma \sqrt{\frac{\log M}{n}} P_n^k(T_n \neq e_k).$$

Thus a sufficient condition for (2.3) to hold is

$$(6.1) \qquad \inf_{T_n} \sup_{k=1,\ldots,M} P_n^k(T_n \neq e_k) \geq c > 0.$$

Since $P_n^k$ is the product of $n$ multivariate Gaussian measures with means $e_k(\sigma/2)\sqrt{\log(M)/n}$ and the covariance matrices $\sigma^2 I$, the Kullback–Leibler divergence between $p_n^k$ and $p_n^1$ is given explicitly by $\mathcal{K}(p_n^k, p_n^1) = (\log M)/4$, for any $k = 2, \ldots, M$, where $p_n^k$ denotes the density of $P_n^k$. We can therefore apply Proposition 2.3 in [25] with $\alpha^* = (\log M)/4$. Taking in that proposition $\tau = 1/M$ we get (6.1) with some $c > 0$ which finishes the proof.

**Acknowledgments.** We would like to thank Jean-Yves Audibert, Arnak Dalalyan and Gilles Stoltz for the remarks that helped to improve the text.

A. JUDITSKY
LABORATOIRE JEAN KUNTZMANN
UNIVERSITÉ GRENOBLE 1
51 RUE DES MATHÉMATIQUES
BP 53
38041 GRENOBLE CEDEX 9
FRANCE
E-MAIL: anatoli.iouditski@imag.fr

P. RIGOLLET
SCHOOL OF MATHEMATICS
GEORGIA INSTITUTE OF TECHNOLOGY
686 CHERRY STREET
ATLANTA, GEORGIA 30332-0160
USA
E-MAIL: rigollet@math.gatech.edu

A. TSYBAKOV
LABORATOIRE DE STATISTIQUE, CREST
3 RUE PIERRE LAROUSSE
92240 MALAKOFF CEDEX
FRANCE
AND
LPMA, CNRS UMR 7599
UNIV. PARIS 6
4 PL. JUSSIEU, CASE 188
75252 PARIS CEDEX 5
FRANCE
E-MAIL: alexandre.tsybakov@ensae.fr